\documentclass[12pt]{amsart}
\usepackage{amsmath, latexsym}
\usepackage[cmtip, arrow]{xy}
\usepackage{amssymb}
\usepackage{lineno}
\usepackage{tikz}
\newtheorem{theorem}{Theorem}[section]
\newtheorem{note}{Note}[section]

\newtheorem{lemma}[theorem]{Lemma}
\newtheorem{corollary}[theorem]{Corollary}

\theoremstyle{definition}
\newtheorem{definition}[theorem]{Definition}
\newtheorem{example}[theorem]{Example}
\newtheorem{remark}[theorem]{Remark}
\newcommand{\R}{\mathbb{R}}                        % real numbers
\newcommand{\C}{\mathbb{C}}                        % complex numbers
\newcommand{\Q}{\mathbb{Q}}                        % rational numbers
\newcommand{\PP}{\mathbb{P}}
\newcommand{\LL}{\mathcal{L}}
\newcommand{\K}{\mathbb{K}}

\def\trdeg{{\rm trdeg}}
\def\trop{{\rm trop}}

\title[Number of connected components of a tropical prevariety]
{Bounds on the number of connected components for
  tropical prevarieties}
\begin{document}
\author[Alex Davydow, Dima Grigoriev]{
\centerline{Alex Davydow}
%\\[-1pt]
\vspace{1mm}
%\small
\centerline{Academic University, ulitsa Khlopina 8,}
\vspace{1mm}
\centerline{Sankt-Petersbourg, 194021, Russia}
% \\[-3pt]
\vspace{1mm}
\centerline{{\tt adavydow@gmail.com}}
%\and
\vspace{3mm}
\centerline{Dima Grigoriev}
% \\[-1pt]
\vspace{1mm}
%\small
\centerline{CNRS, Math\'ematiques, Universit\'e de Lille}
% \\[-3pt]
%\vspace{1mm}
%\small
\vspace{1mm}
\centerline{Villeneuve d'Ascq, 59655, France}
%\\[-3pt]
\vspace{1mm}
\centerline{{\tt \small Dmitry.Grigoryev@math.univ-lille1.fr}}
%\\[-3pt]
\vspace{1mm}
%\small
\centerline{http://en.wikipedia.org/wiki/Dima\_Grigoriev}}
\maketitle
\begin{abstract}
For a tropical prevariety in ${\R}^n$ given by a system of $k$
tropical polynomials in $n$ variables with degrees at most $d$, we
prove that its number of connected components is less than ${k+7n-1
\choose 3n} \cdot \frac{d^{3n}}{k+n+1}$. On a number of $0$-dimensional
connected components a better bound ${k+4n \choose 3n} \cdot
\frac{d^n}{k+n+1}$ is obtained, which extends the Bezout bound due
to B.~Sturmfels from the the case $k=n$ to an arbitrary $k\ge n$.
Also we show that the latter bound is close to sharp, in particular,
the number of connected components can depend on $k$.
\end{abstract}

\tableofcontents

{\bf AMS classification}: 14T05

\bibliographystyle{unsrt}
\section{Introduction}
Let a tropical prevariety $V\subset \R^n$ (see e.g. \cite{richter2005first}) be
given by $k$ tropical polynomials $f_1,\dots,f_k$ in $n$ variables
with (tropical) degrees at most $d$. The principal motivation of
this paper is to bound the number $c$ of connected components of
$V$. Recall (see e.g. \cite{richter2005first}) that $V$ is a polyhedral complex. The
main result (Corollary~\ref{corollary-connected-bound-infinite}) states the bound
\begin{eqnarray}\label{connected}
c\le {k+7n-1 \choose 3n} \cdot \frac{d^{3n}}{k+n+1}
\end{eqnarray}

For the number of isolated points of $V$ (being its $0$-dimensional
connected components) we obtain
(Corollary~\ref{isolated}) a better bound
\begin{eqnarray}\label{bezout}
{k+4n \choose 3n} \cdot \frac{d^n}{k+n+1}
\end{eqnarray}
It can be treated as a generalization of the Bezout inequality on the
number of {\it stable} solutions (see \cite{richter2005first},
\cite{tabera2008tropical} and Section~\ref{section-stability} below)
proved in the case $k=n$ to the case of overdetermined (i.~e. $k>n$)
tropical systems. Recall that $k\ge n$ in order $V$ to have an
isolated point since the local codimension at any point of $V$ is less
or equal to $k$ \cite{shafarevich1977basic}, see also Theorem~\ref{theorem-codimension}. Moreover,
\cite{richter2005first} have proved a tropical Bezout theorem which
states that the number of stable solutions (counted with
multiplicities) of $n$ tropical polynomials $f_1,\dots,f_n$ with
degrees $d_1,\dots,d_n$ respectively, equals $d_1\cdots d_n$.

In Section~\ref{section-lower-bound} we show that bound (\ref{bezout})
is close to sharp by an explicit construction of tropical systems.

The observed phenomenon of dependency of the number of connected
components on $k$ in (\ref{connected}) and in (\ref{bezout}) occurs
similarly for real semialgebraic sets (moreover, for the sum of
Betti numbers which strengthens the bounds established by Oleinik-Petrovskii, Milnor, Thom)
 \cite{basu1996bounding}, while due to a different
reason.

Note that in the case of an algebraic variety given by a polynomial
system $g_1=\dots=g_k=0$ where the degrees of polynomials in $n$
variables do not exceed $d$, the sum of the degrees of the irreducible
components of the variety is bounded by $d^n$, i.~e.  does not depend
on $k$. This holds because the variety remains the same if to replace
$g_1,\dots,g_k$ by their $n+1$ generic linear combinations (see
e.g. \cite{chistov1986algorithm} and \cite{grigoriev1986polynomial}).

Our conjecture is that the sum of Betti numbers of a tropical
prevariety $V$ is bounded by (\ref{connected}). In Theorems~\ref{bet1}, \ref{bet2}
one can find somewhat weaker bounds on the sum of Betti numbers.

The important technical tool to study a system of tropical polynomials
(see Section~\ref{section-vertex}) is the {\it star table} (exploited
in \cite{grigoriev2013complexity}, \cite{grigoriev2012complexity})
consisting of the set of monomials from the given tropical system in
which the minimum is attained at a given point $v\in V$ (here a
monomial is treated as a classical linear function). In these terms we
define a {\it generalized vertex} $v$ of $V$ when the star table is
maximal under inclusion. We produce a description of generalized
vertices in terms of the exponents vectors of the starred monomials
(Theorem~\ref{theorem-vertex-linear-space}).  Then we prove that any
connected component of a tropical prevariety given by a system of tropical
polynomials of fixed degrees with all finite coefficients  contains a generalized vertex
(Theorem~\ref{theorem-connected}).

In Section~\ref{section-stability} we study stable points of a
tropical prevariety given by $n$ tropical polynomials, and provide a
criterion to be a stable point again in terms of the exponent vectors
of the starred monomials
(Theorem~\ref{theorem-stability-linear-space}). This implies that a
generalized vertex of $V$ is a stable point of a suitable {\it
  multisubset} of $\{f_1,\dots,f_k\}$, consisting of $n$ elements
(Theorem~\ref{theorem-vertex-stability}). The established results provide a slightly better
bound than (\ref{connected}) in case of finite coefficients
(Corollary~\ref{corollary-connected-bound}).

To get a bound slghtly better than (\ref{bezout}) in case of finite coefficients
we prove in Section~\ref{section-bezout}
that an isolated point of $V$ is a stable point of an appropriate
subset consisting of $n$ elements among $\{f_1,\dots,f_k\}$
(Theorem~\ref{theorem-choose}). We emphasize that here we consider a
subset, rather than a multisubset as in
Theorem~\ref{theorem-vertex-stability}, this explains the difference
between bounds (\ref{connected}) and (\ref{bezout}).

In Section~\ref{section-compactification} we show that adding $n$ extra variables and $2n$ extra
tropical polynomials to $\{f_1,\dots,f_k\}$ we get a compact
tropical prevariety being homotopy equivalent to $V$. Thus, the
problem of bounding the number of connected and moreover, the sum of Betti numbers of $V$ reduces to a
compact tropical prevariety. Also in Section~\ref{section-compactification} we discuss systems of
tropical polynomials with coefficients allowed to include infinity which allows one to complete the
proofs of bounds (\ref{connected}) and (\ref{bezout}).

\section{Tropical semi-ring and tropical prevarieties}
\begin{definition}
Semi-fields $\R$ and $\R_{\infty}=\R\cup \{\infty\}$ endowed with
operations $\oplus:= \min,\, \otimes:= +,\, \oslash:=-$ are called
tropical semi-rings with or without infinity correspondingly.
\end{definition}

We will denote tropical semi-rings with or without infinity as $\K$
and $\K_{\infty}$ correspondingly.

As well as tropical addition and multiplication we will use
tropical power: $x^{\otimes i}:= x\otimes \cdots \otimes x$.

In this paper we will study tropical polynomials and at first we have
to define a tropical monomial:

\begin{definition}
\emph{Tropical monomial} $Q$ is defined as $Q=a\otimes x_1^{\otimes
  i_1}\otimes \cdots \otimes x_n^{\otimes i_n} = a+i_1\cdot
x_1+\cdots+i_n\cdot x_n$, its \emph{tropical degree}
$\trdeg=i_1+\cdots +i_n$.
\end{definition}

\begin{note}
  As in classic monomials we will often omit multiplication sign when
  it is clear if we speak about a tropical multiplication or a classic
  one. In addition we will omit multiplier of $0$ as it is neutral
  element of tropical multiplication.
\end{note}

Now we can define a tropical polynomial:

\begin{definition}
\emph {Tropical polynomial} $f$ is defined as $f=\bigoplus_j
(a_j\otimes x_1^{\otimes i_{j1}}\otimes \cdots \otimes x_n^{\otimes
  i_{jn}})=\min _j \{Q_j\}$;

$x=(x_1,\dots,x_n)\in \R_{\infty}^n$ is a {\bf tropical zero} of $f$
if at point $x$ either minimum $\min_j \{Q_j(x)\}$ is attained for at
least two different values of $j$ if $\min_j \{Q_j(x)\}$ is finite or
$Q_j(x) = \infty$ for all $j$. If $x\in \R^n$ we say that the tropical
zero is finite. If all monomials with $\trdeg \le d$ present at $f$ we
say that polynomial $f$ of tropical degree $d$ has all its
coefficients finite.
\end{definition}

Then we define a tropical hypersurface:

\begin{definition}
The set of tropical zeros from $\R^n$ of a tropical polynomial is
called a \emph{tropical hypersurface}.
\end{definition}

And finally a tropical prevariety:
\begin{definition}
\emph{Tropical prevariety} is the intersection of a finite number of
tropical hypersurfaces.
\end{definition}

\section{Hahn series and tropical varieties}
To introduce tropical varieties it will be convenient to use a
generalization of Puiseux series known as Hahn (or
Hahn-Mal'cev-Neumann) series (see \cite{hahn1995nichtarchimedischen}).
\begin{definition}
The field of \emph{Hahn series} $K[[T^\Gamma]]$ in the indeterminate
$T$ over a field $K$ and with value group $\Gamma$ (an ordered group)
is the set of formal expressions of the form $f = \sum\limits_{e \in
  \Gamma}{c_eT^e}$ with $c_e \in K$ such that the support $\{e \in
\Gamma : c_e \ne 0\}$ of $f$ is well-ordered. The sum and product of
$f = \sum\limits_{e \in \Gamma}{c_eT^e}$ and $g = \sum\limits_{e \in
  \Gamma}{d_eT^e}$ are given by $f + g = \sum\limits_{e \in
  \Gamma}{(c_e + d_e)T^e}$ and $fg = \sum\limits_{e \in
  \Gamma}{\sum\limits_{e^{'} + e^{''} = e} {c_{e^{'}} d_{e^{''}}
    T^e}}$ ($\sum\limits_{e^{'} + e^{''} = e}$ is finite as
well-ordered set could not contain infinite decreasing sequence).
\end{definition}

To define a tropical variety we have to introduce the operation of
tropicalization.

\begin{definition}
  \emph{Tropicalization} of $x^{'} \in K[[T^\Gamma]]$ is a point $x
  \in \Gamma \cup \{\infty\}$ equal to the least power of $T$ in
  $x^{'}$ if $x^{'}$ is not equal to zero, or $\infty$ otherwise.

  We will denote operation of tropicalization by $\trop$.

  \emph{Tropicalization} $V$ of a variety $V^{'}$ over the field of
  Hahn series $K[[T^\Gamma]]$ consists of the closure in the euclidean topology of the
set of points $x \in \Gamma^n$
  %(\Gamma \cup \{\infty\})^n$
for which there is a point
  $x^{'}=(x^{'}_1, x^{'}_2, \cdots x^{'}_n) \in V^{'}$ with $x^{'}_1 \cdots x^{'}_n \neq 0$, such that
  $x=(\trop(x^{'}_1), \trop(x^{'}_2), \cdots \trop(x^{'}_n)) $.  Set
  $V$ is referred to as a \emph{tropical variety}.
\end{definition}

While any tropical hypersurface is a tropicalization of a hypersurface
over the field of Hahn series $\C[[T^{\R}]]$ (cf. \cite{kapranov}) some tropical
prevarieties do not correspond to any varieties over
$\C[[T^{\R}]]$. For example a tropical prevariety given by the linear
system $$ A = \begin{cases} 0 \oplus x \oplus y \oplus z\\ 0 \oplus x
  \oplus 1y \oplus 1z
      \end{cases}. $$ is not a tropical variety. However, any tropical
variety is a tropical prevariety and moreover the following theorem
holds (see \cite{richter2005first}, \cite{bogart}):

\begin{theorem}
  \label{theorem-prevariety}
  For any variety $V^{'}$ given by polynomial system $A^{'}$ in Hahn
  series $\C[[T^{\R}]]^n$ its tropicalization $V$ is a tropical
  prevariety in $R_{\infty}^n$, and $V$ coincides with the intersection
  of tropical hypersurfaces being tropicalizations of all the
  polynomials from the ideal generated by $A'$. Moreover, $V$ equals
  the tropical prevariety determined by the intersection of a finite
  number of tropicalizations of hypersurfaces provided by polynomials
  from the ideal generated by $A'$ (such a finite subset is called a
  tropical basis of the ideal)
\end{theorem}

To study tropical prevarieties we will use some properties of Hahn
series.

\begin{theorem}[\cite{maclane1939universality}]
  For any algebraically closed field $K$ and ordered divisible group
  $\Gamma$ field of Hahn series $K[[T^{\Gamma}]]$ is algebraically closed.
\end{theorem}

Thus (see e.g. \cite{fulton2012intersection}) we can apply Bezout
theorem to $\C[[T^{\R}]]$.

\begin{definition}
  Let $n$ projective hypersurfaces be given in $\PP^{n}(\C[[T^{\R}]])$
  by $n$ homogeneous polynomials in $n + 1$ variables. Point $x$ is a
  \emph{stable intersection point of these hypersurfaces with
    multiplicity $e$} if under generic small perturbation of the
  coefficients of given polynomials corresponding hypersurfaces will
  have exactly $e$ intersection points in a small neighborhood of $x$.
\end{definition}

\begin{theorem}[Bezout's theorem]
  Let $n$ projective hypersurfaces be given in $\PP^{n}(\C[[T^{\R}]])$
  by $n$ homogeneous polynomials in $n + 1$ variables, of degrees
  $d_1, d_1, \cdots, d_n$. Then the number of stable intersection
  points of these hypersurfaces is equal to $d_1 d_2 \cdots d_n$.
\end{theorem}

Another important property of the field of Hahn series $\C[[T^{\R}]]$
implied by the fact that it is algebraically closed is

\begin{theorem}[Dimension of intersection \cite{shafarevich1977basic}]
   Let variety $V'$ be given by a polynomial system $A$ in $n$
   variables over the field of Hahn series $\C[[T^{\R}]]$. Then if
   system $A$ consists of $k$ polynomials the codimension of each
   irreducible component of $V'$ is less or equal to $k$.
\end{theorem}

This properties of Hahn series are important for studying tropical
varieties and prevarieties due to the following theorem:

\begin{theorem}[\cite{bogart2012obstructions}, \cite{bogart}]
  \label{theorem-codimension}
  For any irreducible variety $V^{'}$ of dimension $m$ over the field
  of Hahn series $\C[[T^{\R}]]$ the local dimension at any point
  $x$ of its tropicalization $V$ is equal to $m$.
\end{theorem}

\begin{remark}
  While Theorem~\ref{theorem-codimension} was known for varieties over
  the field of Puiseux series, the proof can be literally extended to
  Hahn series.
\end{remark}

\section{Generalized vertices}
\label{section-vertex}
To study tropical prevarieties it will be convenient to use the
following definition of vertex:

\begin{definition}
  By a \emph{vertex} of a tropical prevariety we will mean a point for
  which we could not choose a direction in such a way that there is a
  neighborhood of the point where prevariety can be represented as a
  generalized open ended cylinder with axis parallel to the chosen
  direction (a generalized open ended cylinder is a product of an
  arbitrary set of a smaller dimension and a line interval).
\end{definition}

In addition we will need a generalization of this definition, and at
first we have to give a definition of a star table of a tropical system
similar to one introduced in \cite{grigoriev2012complexity}:

\begin{definition}
  Let $A$ be a tropical polynomial system of $k$ polynomials in $n$
  variables with the greatest degree $d$ . We associate with it a
  table $A^{*x}$ of the size $k \times {{n + d - 1} \choose d}$ with
  rows corresponding to polynomials and columns corresponding to all
  possible monomials of degree at most $d$ in $n$ variables. We put $*$ to the
  entry $(i,j)$ iff the $j$-th monomial treated as a (classical)
  linear function attains a minimal value among all the monomials at
  the point $x$ in $i$-th polynomial and we leave all others entries
  empty (see Example~\ref{example-stared-table}).
    \begin{example}
    \label{example-stared-table}
    Consider a tropical system
    $$ A = \begin{cases}
      0 \oplus 1x \oplus y \\
      0 \oplus -2x \oplus -2y \oplus -2x^{\otimes 2} \oplus -3xy \oplus
      -1y^{\otimes 2}
      \end{cases}. $$
    At point $(-1,0)$ this system is equal to
    $$ A = \begin{cases}
      0 \oplus 0 \oplus 0 \\
      0 \oplus -3 \oplus -2 \oplus -4 \oplus -4 \oplus -1
      \end{cases},$$
    so
    $$ A^{*(-1,0)} =
    \left[\begin{smallmatrix}
      0 & x & y & x^{\otimes 2} & xy & y^{\otimes 2} \\
      \hline\\
      *&*&*&&&\\
      &&&*&*&
    \end{smallmatrix}\right].$$
    At point $(1,0)$ this system is equal to
    $$ A = \begin{cases}
      0 \oplus 2 \oplus 0 \\
      0 \oplus -1 \oplus -2 \oplus 0 \oplus -2 \oplus -1
      \end{cases},$$
    so
    $$ A^{*(1,0)} =
    \left[\begin{smallmatrix}
      0 & x & y & x^{\otimes 2} & xy & y^{\otimes 2} \\
      \hline\\
      *&&*&&&\\
      &&*&&*&
    \end{smallmatrix}\right].$$
    At point $(-2,-1)$ this system is equal to
    $$ A = \begin{cases}
      0 \oplus -1 \oplus -1 \\
      0 \oplus -4 \oplus -3 \oplus -6 \oplus -6 \oplus -3
      \end{cases},$$
    so
    $$ A^{*(-2,-1)} =
    \left[\begin{smallmatrix}
      0 & x & y & x^{\otimes 2} & xy & y^{\otimes 2} \\
      \hline\\
      &*&*&&&\\
      &&&*&*&
    \end{smallmatrix}\right].$$
  \end{example}
\end{definition}

All local properties of the tropical prevariety can be expressed in
terms of this table (see the next theorem). In the next chapter we will show how to test
stability of a solution of a tropical system using this table
(Theorem~\ref{theorem-stability-linear-space}), for another example
see \cite{grigoriev2012complexity} where the star table is used to calculate the local
dimension of a linear prevariety.

\begin{theorem}
  \label{theorem-star-table-local}
  Let a tropical prevariety $V$ be given by a tropical system $A$.  If
  $A^{*y} = A^{*z}$ then there is an $\epsilon$, such that
  $\epsilon$-neighborhood of point $y$ of $V$ is homeomorphic to
  $\epsilon$-neighborhood of point $z$ of $V$, moreover this
  homeomorphism is given by a shift of coordinates which sends $y$ to
  $x$.
  \begin{proof}
    Let $d$ be the maximal degree of polynomials in $A$ and $x = (x_1,
    x_2, \cdots. x_n)$ be a set of variables of these polynomials.
    Let's denote by $\Delta_y$ and $\Delta_z$ the minimal differences
    between the values of the starred and non-starred monomials from
    the same polynomials at point $y$ and point $z$
    correspondingly. Let $\Delta=\min(\Delta_y, \Delta_z)$. Denote
    $\epsilon=\frac{\Delta}{3d}$.

    Now we will prove that $\epsilon$ fits the requirements of the
    theorem.  Let's make a change of variables $x'_i=x_i-y_i$ which
    corresponds to a shift of tropical prevariety in such a way that
    $y$ is shifted to $0$. The resulting tropical system we denote by
    $B$. Let's denote the shifted prevariety by $W$. Due to our choice
    of $\epsilon$ in $\epsilon$-neighborhood of $0$ only monomials
    which are starred at $0$ can be starred (they are not greater than
    $0+d\epsilon=\frac{\Delta}{3}$, while others are not lesser than
    $\Delta - d\epsilon=\frac{2\Delta}{3}$), so while studying $B$ in
    $\epsilon$-neighborhood of $0$ we can w.l.o.g. assume that all
    non-starred monomials are infinite. Moreover w.l.o.g. we can
    assume that all coefficients in starred monomials in $B$ are equal
    to zero, otherwise we can tropically multiply corresponding
    polynomials to change them to zero (see
    Example~\ref{example-shift-to-zero}).

    Now if we repeat the same operation replacing all occurrences of
    $y$ by $z$ we will obtain the system which will be the same as $B$
    up to assumptions we made in the end of the previous paragraph. So
    $\epsilon$-neighborhood of $z$ can be obtained from
    $\epsilon$-neighborhood of $y$ by a shift (as both of them can be
    obtained by a shift from $\epsilon$-neighborhood of $0$ of $W$).
  \end{proof}
\end{theorem}

\begin{example}
  \label{example-shift-to-zero}
  Consider system
  $$ \begin{cases}
    2x^{\otimes 2} \oplus x \oplus 0
  \end{cases}. $$
  Assume that we want to study the prevariety given by this system in
  the neighborhood of the point $x = -2$. First we make a change of
  variable $x' = x + 2$:
  $$ \begin{cases}
    -2x'^{\otimes 2} \oplus -2x' \oplus 0
  \end{cases}, $$
  then tropically multiply the polynomial by 2:
  $$ \begin{cases}
    x'^{\otimes 2} \oplus x' \oplus 2
    \end{cases}. $$
  And in $\frac{1}{3}$-neighborhood of $0$ this system can be replaced by:
  $$ \begin{cases}
    x'^{\otimes 2} \oplus x'
    \end{cases}. $$
\end{example}

Now we can give a definition of a generalized vertex:

\begin{definition}
  A point $x$ is a \emph{generalized vertex} of a tropical polynomial
  system $A$ iff the star table $A^{*x}$ is strictly maximal with
  respect to inclusion, i.e. for any other point $y \ne x$ the star
  table $A^{*y}$ does not contain $A^{*x}$.
\end{definition}

\begin{example}
  Point $(-2, -1)$ is not a generalized vertex for a system $A$ from
  Example~\ref{example-stared-table} as $A^{*(-1,0)}$ is greater than
  $A^{*(-2,-1)}$ with respect to inclusion.
\end{example}

\begin{theorem}
  \label{theorem-generalized-vertex-criteria}
  A point $x$ is a \emph{generalized vertex} of a tropical polynomial
  system $A$ iff there is no vector along which the directional
  derivative of every starred monomial in $A^{*x}$ in every polynomial
  is the same (starred monomials from different polynomials can have
  different directional derivatives, see
  Example~\ref{example-nonvertex-movement}), i. e. there is no line
  that passes through the point $x$ along which we can move while
  preserving star table the same in some neighborhood of $x$.
  \begin{proof}
    \begin{enumerate}
    \item First we prove, that if there is a vector along which the
      directional derivative of every starred monomial in $A^{*x}$ in
      every polynomial is the same, then $x$ is not a generalized
      vertex.  It's so because if we will move from point $x$ in the
      direction of this vector there will be a neighborhood where we
      will preserve the star table (so initial star table was not
      strictly maximal).
    \item Now we prove the converse: if $x$ is not a generalized
      vertex then there is a vector along which directional derivative
      of every starred monomial in $A^{*x}$ in every polynomial is the
      same. If $x$ is not a generalized vertex, then there is a point
      $y$ whose star table $A^{*y}$ contains $A^{*x}$. Directional
      derivative along the vector $(x, y)$ will be the same for all
      points stared in $A^{*x}$ in each polynomial, because difference
      between these monomials' values in the same polynomial is the same
      (they are equal both in point $x$ and point $y$).
    \end{enumerate}
  \end{proof}
\end{theorem}

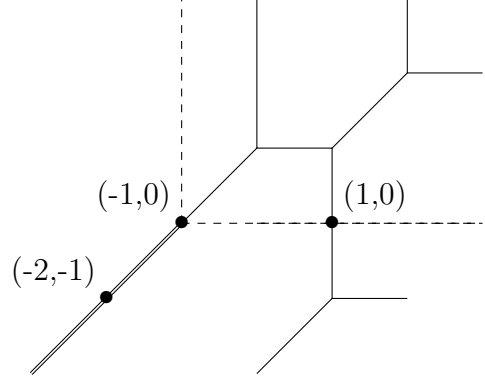
\begin{figure}
  \centering
  \caption {Illustration for Example~\ref{example-nonvertex-movement}}
  \begin{tikzpicture}
    \draw
    (0,1) -- (-1,0)
    (0,1) -- (1,1)
    (1,1) -- (1,-1)
    (1,-1) -- (0,-2)
    (1,-1) -- (2,-1)
    (1,1) -- (2,2)
    (2,2) -- (3,2)
    (2,2) -- (2,3)
    (0,1) -- (0,3);
    \draw[dashed]
    (-1,0) -- (-1,3)
    (-1,0) -- (3,0)
    (0,0) -- (3,0);
    \draw[double]
    (-1,0) -- (-3, -2);
    \draw (-1,0) node {\textbullet};
    \draw (-1,0) node[anchor=south east] {(-1,0)};
    \draw (1,0) node {\textbullet};
    \draw (1,0) node[anchor=south west] {(1,0)};
    \draw (-2,-1) node {\textbullet};
    \draw (-2,-1) node[anchor=south east] {(-2,-1)};

  \end{tikzpicture}
  \label{picture-nonvertex-movement}
\end{figure}

\begin{example}
  \label{example-nonvertex-movement}
  Consider system $A$ from Example~\ref{example-stared-table}.  Point
  $(-2, -1)$ is not a generalized vertex, because we can choose a
  vector $(1,1)$, and directional derivatives of all starred monomials
  in the first polynomial along this vector will be the same and equal
  to $\frac{1}{\sqrt{2}}$, while directional derivatives of starred
  monomials in the second polynomial along this vector will be the same
  and equal to $\sqrt{2}$.

  Points $(1,0)$ and $(-1,0)$ are generalized vertices because we
  could not find a vector with the required property.

  The corresponding prevariety is drawn on
  Figure~\ref{picture-nonvertex-movement}. Prevariety is depicted with
  double lines, the first hypersurface with dashed and the second one
  with solid.
\end{example}

Let's define function $p_n(m_1)$ which takes a tropical monomial in
$n$ variables $x_1,x_2,\cdots,x_n$ as arguments and returns a vector
in $\R^n$ in the following way: $p_n(cx_1^{\otimes a_1}x_2^{\otimes a_2} \cdots
x_n^{\otimes a_n}) = (a_1, a_2, \cdots, a_n)$.

Let's define function $v_n(m_1, m_2)$ which takes two tropical
monomials in $n$ variables as arguments and returns a vector in $\R^n$
in the following way: $v_n(m_1, m_2) = p_n(m_2) - p_n(m_1)$.
\begin{example}
  \begin{itemize}
  \item $v_2(0, x_1) = (1, 0)$,
  \item $v_2(x_1^{\otimes 2}x_2, x_1x_2^{\otimes2}) = (-1, 1)$,
  \item $v_3(0, 2x_1x_2x_3) = (1, 1, 1)$.
  \end{itemize}
\end{example}

Now we can give a criterion of a point being a generalized vertex in
terms of the star table just of this point invoking also functions $v_n$.

\begin{theorem}
  \label{theorem-vertex-linear-space}
  Assume that for a tropical polynomial system $A$ of $k$ polynomials in $n$
  variables with a solution at $x$ we can choose $2n$ monomials
  $m_{i, j}, 1 \le i \le n, 1 \le j \le 2$ with the following properties:
  \begin{itemize}
  \item one monomial can be chosen several times.
  \item the monomials $m_{i,1}$ and $m_{i,2}$ are marked with a
    star in $A^{*x}$ in the same line,
  \item the linear span of vectors $v_n(m_{1,1}, m_{1,2}),
    v_n(m_{2,1},m_{2,2}), \cdots, v_n(m_{n,1}, m_{n,2})$ has dimension
    equal to $n$.
  \end{itemize}
  Then $x$ is a generalized vertex of this system, conversely if $x$
  is a generalized vertex we can always choose $2n$ monomials in the
  described way.
  \begin{proof}
    \begin{enumerate}
      \item First we prove in one direction: if we can not choose $2n$
        monomials with the required properties, then $x$ is not a
        generalized vertex. If we could not choose $2n$ monomials with
        the required properties, then it would mean that the linear
        span of vectors $v_n(y, z)$ where $y$ and $z$ are arbitrary
        starred monomials from the same polynomial has dimension (over
        all the polynomials) lesser than $n$. But if we choose a
        vector orthogonal to this linear span, directional derivatives
        of any pair of starred monomials from the same polynomial
        along this vector will be the same (as directional derivative
        of their difference will be equal to zero). And by
        Theorem~\ref{theorem-generalized-vertex-criteria} this means
        that $x$ is not a generalized vertex.
      \item Now we prove the converse: if $x$ is not a generalized
        vertex then we can not choose $2n$ monomials with the required
        properties. By
        Theorem~\ref{theorem-generalized-vertex-criteria} we can
        choose a vector $v$ along which all directional derivatives of
        starred monomials from the same polynomial will be the same
        (for all the polynomials from $A$). And this means that $v$ is
        orthogonal to $v_n(y, z)$ for any monomials $y$ and $z$
        starred in the same polynomial. But the latter contradicts to that
        the dimension of the linear span of vectors $v_n(m_{1,1},
        m_{1,2}), v_n(m_{2,1},m_{2,2}), \cdots, v_n(m_{n,1}, m_{n,2})$
        equals $n$.
    \end{enumerate}
  \end{proof}
\end{theorem}

Generalized vertex is indeed a generalization of a vertex:

\begin{theorem}
  \label{theorem-vertex-criterion}
  If $x$ is a vertex point of the prevariety $V$ given by a tropical
  polynomial system $A$ then $x$ is a generalized vertex point of $A$.
  \begin{proof}
    Assume the contrary. Then by
    Theorem~\ref{theorem-generalized-vertex-criteria} we can choose a
    vector along which all directional derivatives of monomials
    starred in $A^{*x}$ will be the same. That means that if we choose
    a line passing through the point $x$ and directed by this vector,
    then we can move along it in both directions while keeping star
    table the same in some neighborhood of point $x$. And by
    Theorem~\ref{theorem-star-table-local} in this neighborhood of
    point $x$ prevariety $V$ is a generalized open-ended cylinder.
  \end{proof}
\end{theorem}

However, converse is not true:

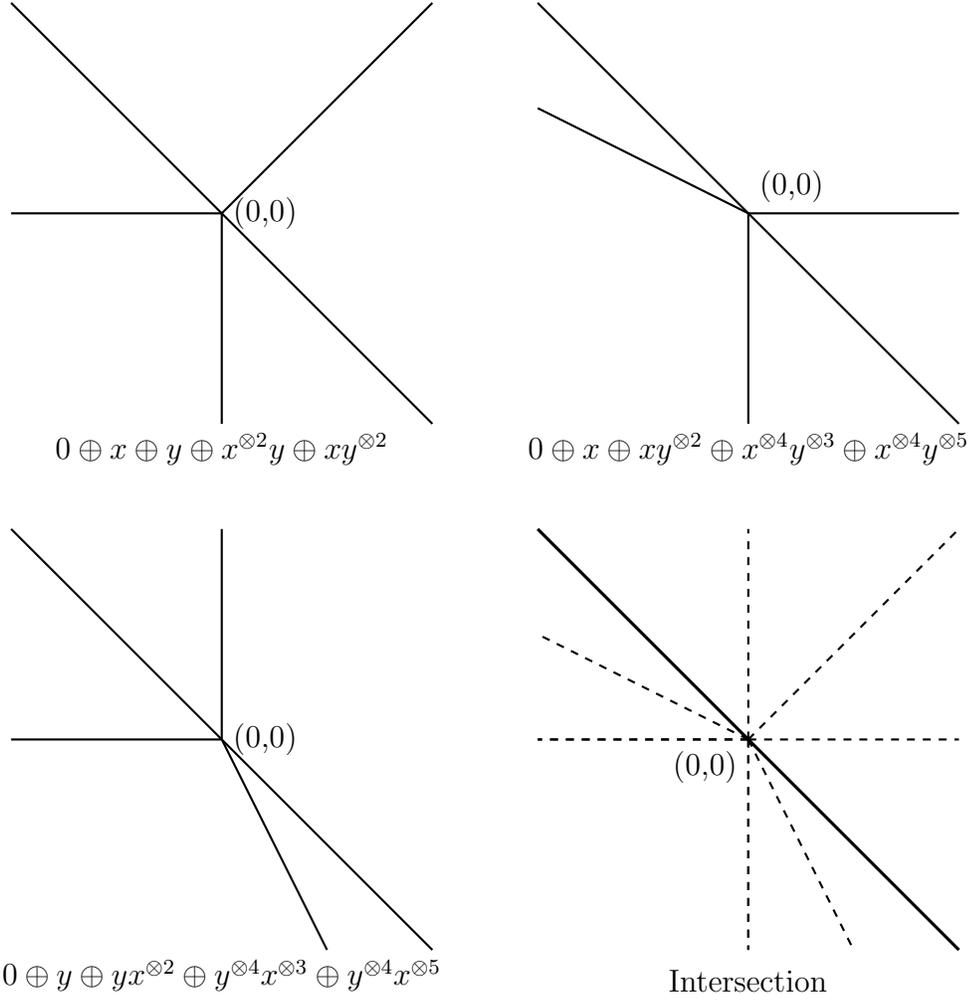
\begin{tikzpicture}[thick, scale=0.7]
\draw[thick] (-5,5) -- ++(4,4);
\draw[thick] (-5,5) -- ++(4,-4);
\draw[thick] (-5,5) -- ++(0,-4);
\draw[thick] (-5,5) -- ++(-4,0);
\draw[thick] (-5,5) -- ++(-4,4);

\draw[thick] (5,5) -- ++(4,0);
\draw[thick] (5,5) -- ++(4,-4);
\draw[thick] (5,5) -- ++(0,-4);
\draw[thick] (5,5) -- ++(-4,2);
\draw[thick] (5,5) -- ++(-4,4);

\draw[thick] (-5,-5) -- ++(0,4);
\draw[thick] (-5,-5) -- ++(4,-4);
\draw[thick] (-5,-5) -- ++(2,-4);
\draw[thick] (-5,-5) -- ++(-4,0);
\draw[thick] (-5,-5) -- ++(-4,4);

\draw[thick, dashed] (5,-5) -- ++(4,4);
\draw[thick, dashed] (5,-5) -- ++(0,-4);
\draw[thick, dashed] (5,-5) -- ++(-4,0);
\draw[thick, dashed] (5,-5) -- ++(4,0);
\draw[thick, dashed] (5,-5) -- ++(0,-4);
\draw[thick, dashed] (5,-5) -- ++(-4,2);
\draw[thick, dashed] (5,-5) -- ++(0,4);
\draw[thick, dashed] (5,-5) -- ++(2,-4);
\draw[thick, dashed] (5,-5) -- ++(-4,0);

\draw[very thick, color=black] (5,-5) -- ++(4,-4);
\draw[very thick, color=black] (5,-5) -- ++(-4,4);

\draw (-5,5) node [anchor=west] {(0,0)};
\draw (5,5) node [anchor=south west] {(0,0)};
\draw (-5,-5) node [anchor=west] {(0,0)};
\draw (5,-5) node [anchor=north east] {(0,0)};

\draw (-5,0) node [anchor=south]
  {$0 \oplus x \oplus y \oplus x^{\otimes 2}y \oplus xy^{\otimes 2}$};
\draw (5,0) node [anchor=south]
  {$0 \oplus x \oplus xy^{\otimes 2} \oplus x^{\otimes 4}y^{\otimes 3}
    \oplus x^{\otimes 4}y^{\otimes 5}$};
\draw (-5,-10) node [anchor=south]
  {$0 \oplus y \oplus yx^{\otimes 2} \oplus y^{\otimes 4}x^{\otimes 3}
    \oplus y^{\otimes 4}x^{\otimes 5}$};
\draw (5,-10) node [anchor=south] {Intersection};

\end{tikzpicture}

\begin{example}
  For a system of tropical polynomials
  $$ \begin{cases}
    0 \oplus x \oplus y \oplus x^{\otimes 2}y \oplus xy^{\otimes 2} \\
    0 \oplus x \oplus xy^{\otimes 2} \oplus x^{\otimes 4}y^{\otimes 3}
    \oplus x^{\otimes 4}y^{\otimes 5} \\
    0 \oplus y \oplus yx^{\otimes 2} \oplus y^{\otimes 4}x^{\otimes 3}
    \oplus y^{\otimes 4}x^{\otimes 5}
    \end{cases}
    $$ $0$ is a generalized vertex, but it is not a vertex (this
  prevariety is equal to a line directed by vector $(-1,1)$).
\end{example}

\section{Stability of Solutions Criteria}
\label{section-stability}
In this paper it will be convenient to use the following definition:
\begin{definition}
  By the \emph{amplitude of a perturbation} we will denote the maximal
  difference between corresponding finite coefficients of the initial
  system and the perturbed one (the infinite coefficients are not
  perturbed).
\end{definition}

In this section we always consider tropical system of $n$ equations in
$n$ variables.

Following Sturmfels, and others \cite{richter2005first} we will use the
following definition of stability:

\begin{definition}
  A point $x$ is a stable point of multiplicity of tropical
  polynomial system $A$ of $n$ equations in $n$ variables if any
  sufficiently small generic perturbation leads to a prevariety being
  a finite set of points, and among them will be exactly $k$ points in
  a neighborhood of $x$.  (More formally: there is $\Delta$, such
  that for any $\delta < \Delta$ there is $\epsilon > 0$ such that for
  any generic perturbation with amplitude not greater than $\epsilon$
  there will be $k$ solutions in $\delta$-neighborhood of $x$).
\end{definition}

This definition can be extended to faces of higher dimension in the
following way:

\begin{definition}
  A $g$-face $L$ is a stable $g$-face of multiplicity $k$ of tropical
  polynomial system $A$ of $n - g$ equations in $n$ variables if any
  sufficiently small generic perturbation leads to a prevariety with
  exactly $k$ $g$-faces in a neighborhood of $x$.  (More formally:
  there is $\Delta$, such that for any $\delta < \Delta$ there is
  $\epsilon > 0$ such that for any generic perturbation with amplitude
  not greater than $\epsilon$ the $\delta$-neighborhood of $x$
  will be intersected by $k$ $g$-faces).
\end{definition}

Our results will be heavily based on the tropical Bezout's equality,
which states the following:
\begin{theorem}[Tropical Bezout's Equality\cite{richter2005first}]
  \label{theorem-bezout-equality}
Every $n$ tropical polynomials with finite coefficients in $n$
variables have $D$ stable finite solutions counted with multiplicities
where $D$ is the product of degrees of the given polynomials.
\end{theorem}

As it was mentioned by Tabera \cite{tabera2008tropical} from this theorem the
following property of stable points of tropical prevariety can be
obtained:

\begin{theorem}
  \label{theorem-stable-perturbation}
  Given $n$ tropical hypersurfaces in $n$-dimensional space the stable
  points of the prevariety being their intersection form a
  well-defined set that varies continuously under perturbations of the
  given hypersurfaces.
\end{theorem}

In this chapter we will always consider
%finite solutions and
systems
with finite coefficients (this means that no monomial of degree at most $d$ in a polynomial
can be omitted), unless we set some of the coefficients to infinity
explicitly.

For effective usage of Theorems~\ref{theorem-bezout-equality}, \ref{theorem-stable-perturbation}
we have to introduce several
simple criteria of stability. While proving theorems we will often
w. l. o. g. study stability at point $0$, and we will consider all
minimal coefficients (i. e. coefficients of the starred monomials in
$A^{*0}$) to be equal to $0$ too (this specific case can be obtained
from the general case by the change of variables to shift point under
consideration to $0$ and by tropical multiplication of equations by
constants, see Example~\ref{example-shift-to-zero}).

\begin{theorem}
  \label{stability-criteria}
  Given a tropical polynomial system $A$ with a solution in $x$, let
  us replace all the coefficients of monomials which are starred in
  $A^{*x}$ by arbitrary set of real numbers and the rest of
  coefficients by infinity (the resulting system denote by $C$). Point
  $x$ is a stable solution of $A$ iff $C$ will have a finite tropical
  solution for any set of chosen real numbers.
  \begin{example}
    Consider tropical system:
    $$ \begin{cases}
      0 \oplus 3x \oplus 0xy \oplus 0x^{\otimes 2} \\
      3 \oplus 0x \oplus 0y^{\otimes 3}
      \end{cases}.$$
    $0$ is a stable solution of this system as system:
    $$ \begin{cases}
      a_1 \oplus a_2xy \oplus a_3x^{\otimes 2} \\
      a_4x \oplus a_5y^{\otimes 3}
      \end{cases}, $$
    has a finite solution for any real $a_1, a_2, a_3, a_4, a_5$.
  \end{example}
  \begin{proof}
    W. l. o. g. we can assume that $x$ is a zero and the coefficients
    of all starred monomials in $A^{*0}$ are equal to $0$.

    Let $d$ be the maximal tropical degree of polynomials in the
    system and let the smallest nonzero coefficient in the initial
    equation be equal to $\Delta$.
    \begin{enumerate}
    \item First we prove in one direction: if $0$ is a stable point of
      $A$, then we can find a finite solution of $C$ for any set of
      coefficients taken as in the theorem. We will prove that for a
      fixed set of coefficients there is a solution. W. l. o. g. we
      can assume that all the coefficients are positive (otherwise we
      can tropically multiply equations by a constant). Let the
      greatest coefficient be equal to $M$. As $0$ is a stable
      solution of the initial system we can choose a $\delta$ with the
      following properties:
        \begin{itemize}
          \item $0 < \delta < \frac{\Delta}{4d}$,
          \item for any perturbation of parameters of the initial
            system with amplitude less or equal to $\delta$ there will
            be a stable solution in a $\frac{\Delta}{4d}$
            neighborhood of $0$.
        \end{itemize}
        Let's consider a perturbation $B$ of the initial system with
        nonzero coefficients unchanged and zero coefficients replaced
        by the corresponding coefficients from $C$ multiplied by
        $\frac{\delta}{M}$.

        By our choice of $\delta$ we can find a solution $y$ of $B$ in
        a $\frac{\Delta}{4d}$ neighborhood of $0$. Monomials which
        were nonzero in $A$ could not be minimal in this solution as
        they are too large. Indeed, as coefficients change is not
        greater than $\delta$ and solution coordinates are less than
        $\frac{\Delta}{4d}$ the value of monomial which was zero in
        $A$ after the perturbation in $y$ is not greater than
        $\delta+d\frac{\Delta}{4d}<\frac{\Delta}{2}$, while the value
        of monomial which was nonzero in $A$, after perturbation in
        $y$ are at least
        $\Delta-\delta-d\frac{\Delta}{4d}>\frac{\Delta}{2}$.

        If we classically multiply the solution and coefficients of the
        equation from $B$ by $\frac{M}{\delta}$, $y$ will be a
        solution for a multiplied system, and if we change all
        coefficients which were nonzero by infinity, solution still
        remains a solution as all the monomials we have changed were
        not minimal. So we have found a solution for $C$.
      \item Now we prove the converse: if we can find a solution of a
        system $C$ for any replacement of the coefficients then $0$ is
        a stable point. We will prove that we can choose such a
        monotone function $p$ that for any perturbation with amplitude
        $\delta<min(p^{-1}(\frac{\Delta}{4d}),\frac{\Delta}{4d})$
        there is a solution in $p(\delta)$ neighborhood of $0$. Let's
        denote the perturbed system by $E$.

        Replace by infinity all monomials in $E$ which are nonzero in
        the initial system $A$. By our assumption this system will
        have a solution, and as it has a solution, it has a solution
        which can be bounded by $2M n! d^n$, where $M$ is a maximal
        coefficient: by Theorem~\ref{theorem-connected} a tropical
        prevariety has at least one generalized vertex and this vertex
        can be obtained as an intersection of $n$ linearly independent
        hyperplanes (classical), given by equality of $n$ pairs of
        tropical monomials taken as in
        Theorem~\ref{theorem-vertex-linear-space}. If we know which
        hyperplanes intersect in this vertex then we can calculate its
        position using Cramer's rule. As all powers in tropical
        equations are integral and do not exceed $d$, as the system of
        equations for the vertex we will obtain a linear system with
        integral coefficients for variables and with a constant part
        which does not exceed $2M$. We can estimate determinant of
        this system as at least $1$ (it's integral and system is not
        degenerate by Theorem~\ref{theorem-vertex-linear-space}) and
        the determinants in the numerators of the Cramer's formula can
        be estimated by $2M n! d^n$. So we can choose $p(\delta) =
        2\delta n! d^n$.

        This solution will be a solution of $E$, as monomials
        corresponding to non-starred monomials of $A$ are too large to
        be minimal (as coefficients change is not greater than
        $\delta$ and solution coordinates are less than $p(\delta)$,
        the value of a monomial which was starred in $A^{*0}$ after
        perturbation in the new solution point is at most
        $\delta+dp(\delta)<\frac{\Delta}{2}$ and the value of a
        non-starred monomial after perturbation in the new solution
        point is at least
        $\Delta-\delta-dp(\delta)>\frac{\Delta}{2}$). For any small
        perturbation we found a solution in a neighborhood of $0$, so
        $0$ is a stable point.
    \end{enumerate}
  \end{proof}
\end{theorem}

Using this theorem we can prove the following lemma:
\begin{lemma}
  \label{lemma-stability-nonozero}
  If $x$ is a stable solution of a tropical system $A$, then for any
  tropical system $F$ and point $y$, if $F^{*y} = A^{*x}$, then $y$ is
  a stable point of $F$.
  \begin{proof}
    W. l. o. g. we can assume that $x$ and $y$ are equal to $0$, and
    that the coefficients of the monomials starred in $A^{*0}$ and
    $F^{*0}$ are equal to zero.

    As $0$ is a stable point of $A$, by
    Theorem~\ref{stability-criteria} if we set all coefficients in the
    monomials which are non-starred in $A^{*0}$ to infinity and
    replace all the coefficients of the starred monomials by arbitrary
    values the obtained system $C$ will have a solution.

    But the result of replacement (system $C$) is the same for systems
    $A$ and $F$, so if in system $F$ we set all coefficients in the
    monomials which are non-starred in $F^{*0}$ to infinity and
    replace all the coefficients of the starred monomials by arbitrary
    values the obtained system will have a solution. And by
    Theorem~\ref{stability-criteria} this means that $0$ is a stable
    solution of $F$.
  \end{proof}
\end{lemma}

This proposition can be strengthened to the following theorem:
\begin{theorem}
  \label{stability-add-zero}
  If $x$ is a stable solution of a tropical system $A$, then for any
  tropical system $F$ and point $y$, if $F^{*y}$ contains $A^{*x}$,
  then $y$ is a stable point of $F$.
  \begin{proof}
    W. l. o. g. we can assume that $x$ and $y$ are equal to $0$, and
    that the coefficients of the monomials starred in $A^{*0}$ and
    $F^{*0}$ are equal to zero.

    Due to the Theorem~\ref{theorem-stable-perturbation} and
    Theorem~\ref{theorem-bezout-equality} we can refer to a stable
    points movement under a parameter perturbation.  We will prove
    that if we change one of nonzero coefficients to zero (the resulting system we denote by $G$),
    still $0$
    remains a stable solution of $G$. The rest will immediately follow from
    Lemma~\ref{lemma-stability-nonozero}.  We will prove by
    contradiction.  Let $0$ be an unstable solution of $G$.
    %the new system.
    Let $\Delta$ be a minimal distance from $0$ to stable
    solutions of $G$.
    %the new system.
    We can choose $\epsilon > 0$ such
    that if we perturb $G$
    %the new system
    with amplitude less than
    $\epsilon$ then every stable solution will move by distance less
    than $\Delta$.

    Now consider a perturbation of $G$
    %the new system
    with new zero
    coefficient replaced by $\epsilon$ and other coefficients unchanged. By
    Lemma~\ref{lemma-stability-nonozero} perturbed system has $0$ as a
    stable solution, but by choice of $\epsilon$ we get a
    contradiction as no stable solution could move to $0$.
  \end{proof}
\end{theorem}

Now we can formulate the last criterion of stability we needed (we
will use functions $v_n$ defined in Section 2):

\begin{theorem}
  \label{theorem-stability-linear-space}
  Assume that for a tropical polynomial system $A$ of $n$ equations in $n$
  variables with a solution at $x$ we can choose $2n$ monomials
  $m_{i, j}, 1 \le i \le n, 1 \le j \le 2$ with the following properties:
  \begin{itemize}
  \item monomials $m_{i,1}$ and $m_{i,2}$ are from $i$-th polynomial
    and they are starred in $A^{*x}$,
  \item linear span of vectors $v_n(m_{1,1},
    m_{1,2}), v_n(m_{2,1},m_{2,2}), \cdots, v_n(m_{n,1}, m_{n,2})$ has
    dimension equal to $n$.
  \end{itemize}
  Then $x$ is a stable solution of system $A$, conversely if $x$ is a
  stable solution we can always choose $2n$ monomials in the way
  described.
\begin{proof}
  W. l. o. g. we can assume that $x$ is equal to $0$, and that the
  coefficients of the monomials starred in $A^{*0}$ are equal to zero.

  \begin{enumerate}
    \item
      First we will prove in one direction: if we could find monomials
      with described properties, then $0$ is a stable point. By
      Theorem~\ref{stability-add-zero} if we prove that $0$ is a
      stable point of a system with all the coefficients of monomials
      except $m_{i, j}, 1 \le i \le n, 1 \le j \le 2$, replaced by say
      $1$ (this system will have only $2$ monomials with zero
      coefficients in each equation), then $0$ is a stable point of
      system $A$. And by Theorem~\ref{stability-criteria} $0$ is
      stable iff the system with coefficients in nonzero monomials
      replaced by infinity will have a solution for any set of
      coefficients replacing zeros. Now we can notice that tropical
      polynomial system with two monomials in each polynomial is just
      a classical linear system (see
      Example~\ref{tropical-classical-example}) and restriction on
      monomials we gave is just a criterion of this system to have
      rank $n$. So as required this system will have a solution for
      any set of coefficients.

    \item Now we will prove the converse: if we could not find
      monomials with described properties, then $0$ is not a stable
      point. We will prove that we can replace all nonzero
      coefficients by infinity and zero coefficients by arbitrary real
      numbers in such a way that the obtained system will have no
      solution and thus, by Theorem~\ref{stability-criteria} $0$ is
      not a stable point.

      Let's replace all zero coefficients by arbitrary real numbers
      which are linear independent over $\Q$ and nonzero coefficients
      by infinity. Consider a solution of this system. Let's choose
      $f_{i, j}, 1 \le j \le 2$ as pairs of starred monomials from
      $i$-th equation (if there are more than two starred monomials,
      we will choose just two arbitrary among them). Linear span of
      $v_n(f_{i,1}, f_{i,2}), 1 \le i \le n$ has dimension lesser than
      $n$ by assumption, so the system of classical linear equations,
      expressing that $f_{i,1}=f_{i,2}, 1 \le i \le n$ will have the
      rank lesser than $n$. This system has rational coefficients of
      the variables, while free terms from its' equations are linear
      independent over $\Q$, so it has no solutions, as otherwise a
      rational linear dependency between these constants could be
      found. So we come to a contradiction and this means that $0$ is
      not a stable point.
  \end{enumerate}

\end{proof}
\end{theorem}

\begin{example}
  \label{tropical-classical-example}
  Tropical polynomial system:
  $$ \begin{cases}
    x_1x_2 \oplus x_1^{\otimes 3} \\
    6x_1^{\otimes 5} \oplus 4x_2^{\otimes 2}
  \end{cases}$$
  is equivalent to classical linear system:
  $$ \begin{cases}
    x_1 + x_2 = 3x_1 \\
    6 + 5x_1 = 4 + 2x_2
  \end{cases} $$
\end{example}

The criterion from Theorem~\ref{theorem-stability-linear-space} of a
point being a stable solution of a tropical system will be used
further in our paper. In fact, this criterion can be tested in
polynomial time, by means of an algorithm which produces a maximal rank
subset of an intersection of two matroids, see
e.g. \cite{schrijver2003combinatorial}, \cite{barvinok1995new}.
\section{Estimating the number of connected components for tropical
  systems with finite coefficients}
\label{section-connected}
Using the theorems from the Section~\ref{section-stability} we can
bound the number of connected components of a tropical prevariety.

As in the previous section we assume that all the coefficients in a
tropical system are finite.

At first we will show that every connected component contains at least
one generalized vertex.

\begin{theorem}
  \label{theorem-connected}
  If a tropical prevariety $V$ is given by a tropical polynomial
  system $A$ with finite coefficients, then in any connected component
  of $V$ there is at least one generalized vertex.
  \begin{proof}
    Consider a point $x$ of $V$. If it is not a generalized vertex
    then by Theorem~\ref{theorem-generalized-vertex-criteria} there is
    a vector along which all the directional derivatives of starred
    monomials in each polynomial are the same. Let's look at the star
    table while moving from point $x$ forward and backward along this
    vector. In some neighborhood of $x$ the star table will not
    change, but at some point a new star has to appear (as the
    coefficients of $A$ are finite all monomials are present, so there
    will be at least one non-starred monomial whose derivative along
    the chosen vector differs from the derivatives of the starred
    monomials in the same polynomial, as there is no vector along
    which the derivatives of all the monomials are the same). Let's
    choose this point as a new $x$. By this procedure we have
    increased the number of stars in $A^{*x}$. Now we can repeat the
    described process. But as there is a finite number of cells in the
    star table, we can't repeat this process up to infinity, so at
    some step the chosen point $x$ must be a generalized vertex.
  \end{proof}
\end{theorem}

To estimate the number of generalized vertices we will prove the
following theorem:

\begin{theorem}
  For any generalized vertex $x$ of a tropical polynomial system in $n$
  variables we can choose a multiset of $n$ polynomials from this
  system in such a way that $x$ is a stable solution for a tropical
  system given by the chosen polynomials (one polynomial can be chosen
  several times). Moreover, if there were $k \ge n$ polynomials in the
  initial system, then we can choose at least $k - n + 1$ different
  multisets of $n$ polynomials with the described properties.
  \label{theorem-vertex-stability}
  \begin{proof}
    \begin{enumerate}
      Existence of one multiset with the described properties
      immediately follows from
      Theorem~\ref{theorem-vertex-linear-space} and
      Theorem~\ref{theorem-stability-linear-space}.

      The second part of the theorem can be proved in the following
      way: consider a multiset $\{p_1, p_2, \cdots, p_n\}$, of $n$
      polynomials with the described properties. As $x$ is a stable
      point of these polynomials we can choose a set of monomials
      $m_{1..n,1..2}$ as described in
      Theorem~\ref{theorem-stability-linear-space}. Now we will with
      each polynomial associate one vector:
      $a_i=v_n(m_{i,1},m_{i,2})$. By
      Theorem~\ref{theorem-stability-linear-space} this vectors form a
      basis in $\R^n$. Now we will prove that for any polynomial
      $p_{n+1}$ which is not represented in the multiset
      $\{p_1,p_2,\cdots,p_n\}$ we can choose a polynomial $p_i$ in
      such a way that $x$ will be a stable point of the system given
      by polynomials
      $\{p_1,p_2,\cdots,p_{i-1},p_{i+1},\cdots,p_n,p_{n+1}\}$. As $x$
      is a solution of $p_{n+1}$ there are at least two monomials
      which are starred in $p_{n+1}$ at the point $x$. Let's denote
      them by $m_{n+1,1}$ and $m_{n+1,2}$ (if there are more than two
      starred monomials in $p_{n+1}$ at the point $x$ we will choose
      an arbitrary pair of starred monomials). As $\{a_1, a_2, \cdots,
      a_n\}$ is a basis, there should be a linear combination of this
      vectors which will be equal to $v_n(m_{n+1,1},m_{n+1,2})$, this
      means that
      $v_n(m_{n+1,1},m_{n+1,2})=ca_i+L(a_1,a_2,\cdots,a_{i-1},a_{i+1},\cdots,a_n)$
      for $c\ne 0$ and a certain $1 \le i \le n$, where $L$ is a
      linear function. So
      $\{a_1,a_2,\cdots,a_{i-1},a_{i+1},\cdots,a_n,v_n(m_{n+1,1},m_{n+1,2})\}$
      is a basis, and by Theorem~\ref{theorem-stability-linear-space},
      this means that $x$ is a stable point of
      $\{p_1,p_2,\cdots,p_{i-1},p_{i+1},\cdots,p_n,p_{n+1}\}$. As we
      can choose at least $k - n$ polynomials which are not included
      in the initial multiset, there are at least $k - n + 1$
      multisets with the properties required in the theorem.
    \end{enumerate}
  \end{proof}
\end{theorem}

As a consequence from this theorem we can obtain:

\begin{theorem}
  The number of generalized vertices of a tropical prevariety in $\R^n$ given by
  $k$ polynomials with tropical degrees bounded by $d$ and with finite coefficients is not greater
  than $\frac{d^n}{k - n + 1} {{k + n - 1} \choose n}$.
  \begin{proof}
    We can choose up to ${{k + n - 1} \choose n}$ different multisets
    of equations and by Theorem~\ref{theorem-bezout-equality} system
    formed by each of these multisets will have at most $d^n$ stable
    points. Moreover each point will be calculated at least $k-n+1$
    times due to Theorem~\ref{theorem-vertex-stability}. This implies
    the required bound.
  \end{proof}
\end{theorem}

By Theorem~\ref{theorem-connected} we can obtain:
\begin{corollary}
  \label{corollary-connected-bound}
  The number of connected components of tropical prevariety in $\R^n$
  given by $k$ polynomials with tropical degrees bounded by $d$ and with finite coefficients is not
  greater than $\frac{d^n}{k - n + 1} {{k + n - 1} \choose n}$.
\end{corollary}

However this bound is not sharp, and while it's rather precise for
considerably overdetermined system (in
Theorem~\ref{theorem-lower-bound} we will show that for
overdetermined systems a bound $\frac{d^n}{k - n + 1}{k\choose n}$ can
be achieved), for underdetermined systems a better bound can be proved:
\begin{theorem}
  The number of connected components of a tropical prevariety in
  $\R^n$ given by $k$ polynomials with tropical degree bounded by $d$
  is not greater than ${{d + n} \choose d}^{2k}$.
  \begin{proof}
    A tropical prevariety is a union of at most ${{d + n} \choose
      d}^{2k}$ convex polyhedra, each of them given by a star table with
    exactly two stars in every row.
  \end{proof}
\end{theorem}

While this bound is not interesting for overdetermined system, for
small $k$ and $d$ comparatively to $n$ it can be much better than the
bound from Corollary~\ref{corollary-connected-bound}.

Now relying on Corollary~\ref{corollary-connected-bound} we can
obtain a bound on the sum of Betti numbers (discrete Morse's theory
states that in compact tropical prevariety it is bounded by the
number of faces of all dimension, see e. g. \cite{forman2002user}).

\begin{theorem}\label{bet1}
  For any $0\le l\le n$ the §l§-th  Betti number of a compact tropical prevariety given by a
  system of $k$ polynomials of maximal degree at most $d$ in $n$ variables does not exceed
  $(\frac{d^n}{k - n + 1} {{k + n - 1} \choose n})^{l+1}$.
  \begin{proof}
    This result immediately follows from the fact that any
    $l$-dimensional face of a compact tropical prevariety contains at
    least $l + 1$ vertices.
  \end{proof}
\end{theorem}

However this result is far from sharp, for example for small $d$ it
can be improved by the bound on the number of faces for \emph{arrangements}
(arrangement is a union of several hypersurfaces, see
e.g. \cite{fukuda1991combinatorial})
\begin{theorem}\label{bet2}
  The sum of Betti numbers of a compact tropical prevariety given by a
  system of $k$ polynomials of maximal degree $d$ in $n$ variables does
  not exceed $3^n + 2^n {k {n+d \choose n}^2 \choose n} + o((k {n+d \choose n}^2)^n)$.
  \begin{proof}
    The number of all faces in the arrangement could be estimated as
    $3^n + 2^n {m \choose n} + o(m^n)$ where $n$ is a dimension
    and $m$ is the number of hypersurfaces (see e.g. Buck's formula in
    \cite{fukuda1991combinatorial}).  Faces of tropical prevariety is a subset of faces
    of arrangement of hypersurfaces, where for every pair of monomials
    from the same polynomial we add a hypersurface where they are equal.
    Thus we obtain the required bound (the number of monomials in each
    polynomial does not exceed ${n+d \choose n}$).
  \end{proof}
\end{theorem}

\section{Tropical Bezout Inequality for Overdetermined Systems}
\label{section-bezout}
While the bound on the number of connected components obtained in
Corollary~\ref{corollary-connected-bound} can be used as a bound on the
number of isolated points, in this particular case it can be slightly
improved. Throughout of this section we consider an algebraic variety
over Hahn series
$V' \in (\C[[T^{\R}]]\setminus{0})^n$

Now we can prove, that:
\begin{theorem}
  \label{theorem-choose}
  Given an overdetermined tropical polynomial system $A$ of $k \ge n$
  equations in $n$ variables with an isolated solution at $x$ we can
  always choose $2k$ monomials $m_{i,j}, 1 \le i \le k, 1 \le j \le 2$
  with the following properties:
  \begin{itemize}
  \item monomials $m_{i,1}$ and $m_{i,2}$ are taken from $i$-th
    polynomial and starred in $A^{*x}$.
  \item the linear span of vectors $v_n(m_{1,1}, m_{1,2}),
    v_n(m_{2,1},m_{2,2}), \cdots, v_n(m_{k,1}, m_{k,2})$ has dimension
    equal to $n$.
  \end{itemize}
  \begin{proof}
    W. l. o. g. we can assume that $x = 0$ (we can always shift a
    prevariety in such a way that it is). In the proof we will refer
    to $v_n(x, y)$ as a vector given by the pair $(x, y)$.

    Let's prove by contradiction. Consider $2k$ monomials $m_{i,j}, 1
    \le i \le k, 1 \le j \le 2$ and number $l$ with the following
    properties:
    \begin{itemize}
    \item monomials $m_{i,1}$ and $m_{i,2}$ are taken from $i$-th
      polynomial and starred in $A^{*0}$,
    \item the linear span of vectors $v_n(m_{1,1}, m_{1,2}),
      v_n(m_{2,1},m_{2,2}), \cdots, v_n(m_{k,1}, m_{k,2})$ has
      dimension equal to $l$,
    \item for any other $2k$ monomials $m'_{i,j}, 1 \le i \le k, 1
      \le j \le 2$, where $m'_{i,1}, m'_{i,2}$ are taken from
      $i$-th polynomial and starred in $A^{*0}$, the linear span of
      vectors $v_n(m'_{1,1}, m'_{1,2}),
      v_n(m'_{2,1},m'_{2,2}), \cdots, v_n(m'_{k,1},
      m'_{k,2})$ has dimension equal  or less than $l$.
    \end{itemize}
    Let's denote the linear span of $v_n(m_{1,1}, m_{1,2}),
    v_n(m_{2,1},m_{2,2}), \cdots, v_n(m_{k,1}, m_{k,2})$ by $\LL$.
    W. l. o. g. let's assume, that $\left\{v_n(m_{1,1}, m_{1,2}),
    v_n(m_{2,1}, m_{2,2}), \cdots, v_n(m_{l,1}, m_{l,2})\right\}$ is a
    basis of $\LL$. Let's notice that if $i > l$ and $m'_{i,1}$ and
    $m'_{i,2}$ are taken from $i$-th polynomial and starred in
    $A^{*0}$ then $v_n(m'_{i,1}, m'_{i, 2})$ is contained in
    $\LL$ (otherwise the linear span of $v_n(m_{1,1}, m_{1,2}),
    v_n(m_{2,1},m_{2,2}), \cdots, v_n(m_{l,1}, m_{l,2}),
    v_n(m'_{i,1},m'_{i,2})$ will be equal to $l + 1$ which
    contradicts our assumption).

    As $v_n(m'_{i,1}, m'_{i, 2})$ is contained in $\LL$ we can express
    it as a linear combination of the basis vectors. Let's notice that
    if $v_n(m_{j,1}, m_{j,2})$ occurs in the expression of
    $v_n(m'_{i,1}, m'_{i, 2})$ with a non-zero coefficient we can swap
    $j$-th and $i$-th polynomials and at the same time swap $m_{j,1},
    m_{j,2}$ and $m'_{i,1}, m'_{i,2}$ and the linear span of
    $v_n(m_{1,1}, m_{1,2}), v_n(m_{2,1},m_{2,2}), \cdots, v_n(m_{k,1},
    m_{k,2})$ will remain the same. So we can choose a set $S$ of $g
    \le l$ polynomials from the first $l$ polynomials in such a way
    that a vector given by any pair of monomials taken from
    polynomials from $S$ and starred in $A^{*0}$ or polynomials
    standing after $l$-th in $A$ will lay in the linear span of the
    vectors $v_n(m_{i,1}, m_{i,2}), i \in S$.

    Let's replace all the monomials, except $m_{i,j}, i \in S, 1 \le j
    \le 2$ by infinity and consider system $B$ of the first $l$
    polynomials (some of them changed at the previous step). Obviously it
    has a solution at $0$. And locally any solution of this system is
    a solution of the initial system (every vector given by a pair of
    monomials starred in $A^{*0}$ from the polynomials from $S$ in any
    polynomial can be represented as a linear combination of vectors
    $v_n(m_{i,1}, m_{i,2}), i \in S$). So the local dimension of the
    prevariety given by system $B$ at $0$ is zero. But it contradicts
    Theorem~\ref{theorem-codimension} as $l < n$.
  \end{proof}
\end{theorem}

We note that the proof of the latter theorem is somewhat similar to the proof of
Theorem 41.1 on the intersection of matroids in \cite{schrijver2003combinatorial}.

Which leads to:

\begin{theorem}
  \label{theorem-bezout-overdetermined}
  The number of isolated solutions of an overdetermined tropical
  polynomial system of $k \ge n$ polynomials in $n$ variables is not greater
  than $\frac{{k \choose n}} {(k - n + 1)} D$, where $D$ is the product of $n$
  greatest degrees of the given polynomials.
\begin{proof}
  By theorems \ref{theorem-choose}, \ref{stability-criteria} we can
  state that any isolated solution is a stable solution of some
  subsystem of the size $n$ (by a sutable shift of variables we can
  always shift the solution to $0$). There are less or equal than ${k
    \choose n}$ subsystems and by Bezout's equality each of them has
  at most $D$ stable solutions. Moreover each solution is counted at
  least $k - n + 1$ times (the reasoning is the same as in
  Theorem~\ref{theorem-vertex-stability}).
\end{proof}
\end{theorem}

Observe that (similar to the end of Section~\ref{section-stability}) the algorithm which
produces a maximal rank subset of intersection of two matroids
allows one to test with polynomial complexity, whether a
given solution of a tropical polynomial system is isolated.

As we will show in the next section, this bound is close to sharp.
\section{Lower Bounds on the Number of Isolated Tropical Solutions}
\label{section-lower-bound}
In this section we will build an example, which shows that Bezout
inequality in case of tropical polynomial systems is close to
sharp. While we will omit some monomials (i.e. we will use infinite
coefficients), example like this can be build with finite coefficients
only (if infinite coefficients are replaced by sufficiently large
finite numbers).
\begin{theorem}
  \label{theorem-lower-bound}
  Given $n$ we can build a series of tropical systems of $k(n-1), k
  \ge 3$ equations in $n \ge 2$ variables of degree $4d, d \ge 1$ in
  such a way that the number of solutions of systems from this set is
  $2(k - 1)^{n - 1} d^n$.
\begin{proof}
  \begin{figure}
    \centering
    \caption {Hypersurface (curve) $H_1$ given by tropical polynomial
      $3 \oplus 1x_1 \oplus x_1x_2 \oplus 1x_2 \oplus x_1x_2^{\otimes 2} \oplus
      2x_2^{\otimes 2}$ and its Newton's polygon. $\alpha$ is equal to $2$ and
      $\beta$ is equal to $3$ in the picture.}

    \begin{tikzpicture}[scale=1.0]
      \draw
      (1,2) -- (0,2)
      (1,2) -- (1,3)
      (1,2) -- (2,1)

      (2,1) -- (3,1)
      (2,1) -- (2,0)

      (1,3) -- (0,3)
      (1,3) -- (2,4)

      (2,4) -- (2,5)
      (2,4) -- (3,4);
      \draw (2,1) node {\textbullet};
      \draw (2,1) node[anchor=south west] {($\alpha$,$\beta$)};
      \draw (2,4) node {\textbullet};
      \draw (2,4) node[anchor=south west] {($\alpha$,$\beta + 3$)};
      \draw
      (5,0) -- (5,2)
      (5,0) -- (6,0)
      (6,0) -- (6,2)
      (5,2) -- (6,2)

      (5,0) -- (6,1)
      (5,1) -- (6,1)
      (5,2) -- (6,1);
    \end{tikzpicture}
    \label{picture-one-equation}
  \end{figure}
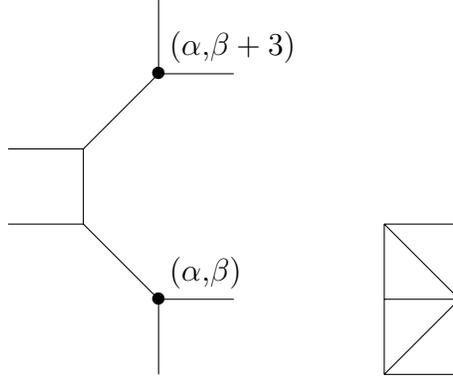

  Consider a tropical polynomial system in 2 variables:
  $$
    A=
    \begin{cases}
      3 \oplus 1x_1 \oplus x_1x_2 \oplus 1x_1 \oplus x_1x_2^{\otimes
        2} \oplus 2x_2^{\otimes 2} \\
      3 \oplus 1x_1 \oplus 3x_1x_2 \oplus 4x_2 \oplus 6x_1x_2^{\otimes
        2} \oplus 8x_2^{\otimes 2} \\
      3 \oplus 1x_1 \oplus 6x_1x_2 \oplus 7x_2 \oplus
      12x_1x_2^{\otimes 2} \oplus 14x_2^{\otimes 2}\\
      \cdots\\
      3 \oplus 1x_1 \oplus (3k - 3)x_1x_2 \oplus (3k -2)x_2 \oplus (6k
      - 6)x_1x_2^{\otimes 2}
        \oplus (6k - 4)2x_2^{\otimes 2}
    \end{cases}
  $$

  The graph of the hypersurface $H_1$ given by the first polynomial is
  depicted on Figure~\ref{picture-one-equation}. The prevariety
  (curve) $H_i$ of the $i$-th polynomial of $A$ is obtained from $H_1$
  by a vertical shift down by $3i-3$.

  Therefore, the points ($\alpha$, $\beta - 3j$), $0 \le j \le k - 2$
  are solutions of $A$.

  Moreover, these points are isolated solutions since the prevariety
  of $A$ consists of these points and of two vertical half-lines.

  \begin{figure}
    \centering
    \caption {Newton's polygon used in
      Theorem~\ref{theorem-lower-bound}}
    \begin{tikzpicture}[scale=1]
      \draw
      (1,0) -- (1,2)
      (1,0) -- (2,0)
      (2,0) -- (2,2)
      (1,2) -- (2,2)

      (1,0) -- (2,1)
      (1,1) -- (2,1)
      (1,2) -- (2,1);

      \draw
      (3,0) -- (3,2)
      (3,0) -- (2,0)
      (2,0) -- (2,2)
      (3,2) -- (2,2)

      (3,0) -- (2,1)
      (3,1) -- (2,1)
      (3,2) -- (2,1);
      \label{picture-square}
    \end{tikzpicture}
  \end{figure}
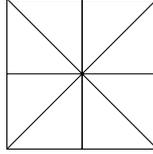

  Now we construct a tropical system $B$ in $2$ variables consisting
  of $k$ polynomials of degrees $4d$ for any $d \ge 1$. The Newton's
  polygon of each of these polynomials is a square with the mesh $2d$
  which is obtained from the $2 \times 2$ square depicted in
  Figure~\ref{picture-square} by replicating it. The coefficients of
  these polynomial are chosen with suitable conditions imposed on the
  distances which follow. The curve of the first polynomial of $B$ is
  depicted on Figure~\ref{picture-many-equations}. The curve consists
  of $d$ horizontal layers of $d$ hexagons each of a height and a
  width equal to $3$ each obtained from the previous one by a vertical
  shift. We impose the condition that the first shift (which is equal
  $\delta - \kappa$) is greater than $3k$. In a similar way the second
  shift $(\kappa - 3) - \nu$ is also greater than $3k$ and so on.  The other
  polynomials are chosen in the way similar to system $A$: they give
  curves which are vertical shifts of the curve given by the first
  polynomial. The second curve is shifted down by $3$, the third curve
  is shifted down by $6$, $\cdots$, the $k$-th curve is shifted down
  by $3(k-1)$.

  Solutions of $B$ form $d$ series of isolated points, each series
  consist of $2(k - 1)d$ points in each series and $2d$ half-lines (for each
  $0\le i\le d-1$) a series has isolated points: $(\gamma + 4i,
  \delta)$, $(\gamma + 4i + 1, \delta)$, $(\gamma + 4i, \delta - 3)$,
  $(\gamma + 4i + 1, \delta - 3)$, $\cdots$, $(\gamma + 4i, \delta -
  3k + 3)$, $(\gamma + 4i + 1, \delta - 3k + 3)$; $(\gamma + 4i,
  (\kappa - 3) - 3)$, $(\gamma + 4i + 1, (\kappa - 3) - 3)$, $(\gamma
  + 4i, (\kappa - 3) - 6)$, $(\gamma + 4i + 1, (\kappa - 3) - 6)$,
  $\cdots$, $(\gamma + 4i, (\kappa - 3) - 3k + 3)$, $(\gamma + 4i + 1,
  (\kappa - 3) - 3k + 3)$; and so on).

  \begin{figure}
    \centering
    \caption {Example for Theorem~\ref{theorem-lower-bound}}

    \begin{tikzpicture}[scale=0.7]
      \draw
      (1,2) -- (0,2)
      (1,2) -- (1,3)
      (1,2) -- (2,1)

      (2,1) -- (3,1);
      \draw [dashed] (2,1) -- (2,0);
      \draw
      (1,3) -- (0,3)
      (1,3) -- (2,4)

      (2,4) -- (3,4)

      (3,4) -- (4,3)

      (4,3) -- (4,2)

      (4,2) -- (3,1);

      \draw [dashed] (3,1) -- (3,0);
      \draw

      (5,2) -- (4,2)
      (5,2) -- (5,3)
      (5,2) -- (6,1)

      (6,1) -- (7,1);
      \draw [dashed] (6,1) -- (6,0);
      \draw
      (5,3) -- (4,3)
      (5,3) -- (6,4)

      (6,4) -- (7,4)

      (7,4) -- (8,3)

      (8,3) -- (8,2)
      (8,3) -- (9,3)

      (8,2) -- (7,1)
      (8,2) -- (9,2);

      \draw [dashed] (7,1) -- (7,0);
      \draw
      (1,6) -- (0,6)
      (1,6) -- (1,7)
      (1,6) -- (2,5)

      (2,5) -- (3,5);
      \draw [dashed] (2,5) -- (2,4);
      \draw

      (1,7) -- (0,7)
      (1,7) -- (2,8)

      (2,8) -- (2,9)
      (2,8) -- (3,8)

      (3,8) -- (3,9)
      (3,8) -- (4,7)

      (4,7) -- (4,6)

      (4,6) -- (3,5);

      \draw [dashed] (3,5) -- (3,4);
      \draw

      (5,6) -- (4,6)
      (5,6) -- (5,7)
      (5,6) -- (6,5)

      (6,5) -- (7,5);
      \draw [dashed] (6,5) -- (6,4);
      \draw

      (5,7) -- (4,7)
      (5,7) -- (6,8)

      (6,8) -- (6,9)
      (6,8) -- (7,8)

      (7,8) -- (7,9)
      (7,8) -- (8,7)

      (8,7) -- (8,6)
      (8,7) -- (9,7)

      (8,6) -- (7,5)
      (8,6) -- (9,6);

      \draw [dashed] (7,5) -- (7,4);

      \draw
      (11,0) -- (11,2)
      (11,0) -- (12,0)
      (12,0) -- (12,2)
      (11,2) -- (12,2)

      (11,0) -- (12,1)
      (11,1) -- (12,1)
      (11,2) -- (12,1);

      \draw
      (11,2) -- (11,4)
      (11,2) -- (12,2)
      (12,2) -- (12,4)
      (11,4) -- (12,4)

      (11,2) -- (12,3)
      (11,3) -- (12,3)
      (11,4) -- (12,3);

      \draw
      (13,0) -- (13,2)
      (13,0) -- (14,0)
      (14,0) -- (14,2)
      (13,2) -- (14,2)

      (13,0) -- (14,1)
      (13,1) -- (14,1)
      (13,2) -- (14,1);

      \draw
      (13,2) -- (13,4)
      (13,2) -- (14,2)
      (14,2) -- (14,4)
      (13,4) -- (14,4)

      (13,2) -- (14,3)
      (13,3) -- (14,3)
      (13,4) -- (14,3);

      \draw
      (13,0) -- (15,2)
      (13,0) -- (12,0)
      (12,0) -- (12,2)
      (13,2) -- (12,2)

      (13,0) -- (12,1)
      (13,1) -- (12,1)
      (13,2) -- (12,1);

      \draw
      (13,2) -- (15,4)
      (13,2) -- (12,2)
      (12,2) -- (12,4)
      (13,4) -- (12,4)

      (13,2) -- (12,3)
      (13,3) -- (12,3)
      (13,4) -- (12,3);

      \draw
      (15,0) -- (15,2)
      (15,0) -- (14,0)
      (14,0) -- (14,2)
      (15,2) -- (14,2)

      (15,0) -- (14,1)
      (15,1) -- (14,1)
      (15,2) -- (14,1);

      \draw
      (15,2) -- (15,4)
      (15,2) -- (14,2)
      (14,2) -- (14,4)
      (15,4) -- (14,4)

      (15,2) -- (14,3)
      (15,3) -- (14,3)
      (15,4) -- (14,3);

      \draw (2,0) node {\textbullet};
      \draw (2,0) node[anchor=east] {($\gamma$,$\nu$)};

      \draw (2,1) node {\textbullet};
      \draw (2,1) node[anchor=east] {($\gamma$,$\kappa - 3$)};

      \draw (2,4) node {\textbullet};
      \draw (2,4) node[anchor=east] {($\gamma$,$\kappa$)};

      \draw (2,5) node {\textbullet};
      \draw (2,5) node[anchor=east] {($\gamma$,$\delta$)};

      \draw (7,0) node {\textbullet};
      \draw (7,0) node[anchor=west] {($\gamma + 5$,$\nu$)};

      \draw (7,1) node {\textbullet};
      \draw (7,1) node[anchor=west] {($\gamma + 5$,$\kappa - 3$)};

      \draw (7,4) node {\textbullet};
      \draw (7,4) node[anchor=west] {($\gamma + 5$,$\kappa$)};

      \draw (7,5) node {\textbullet};
      \draw (7,5) node[anchor=west] {($\gamma + 5$,$\delta$)};

    \end{tikzpicture}
    \label{picture-many-equations}
  \end{figure}
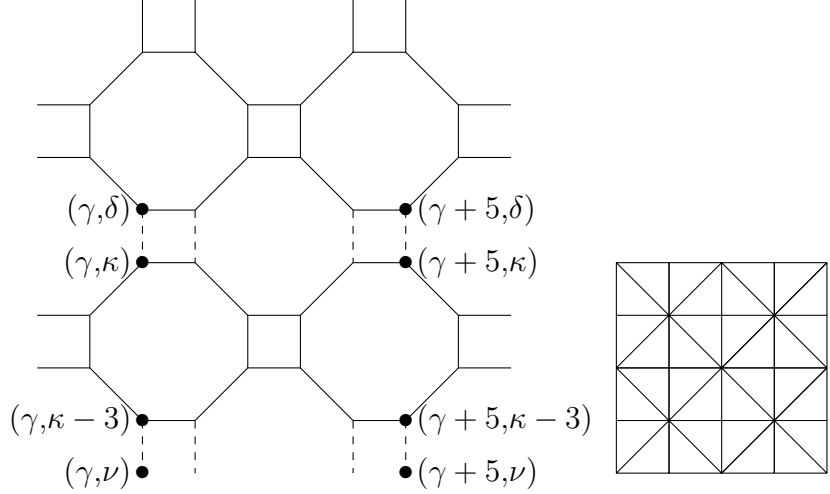
  Now consider system $C_n$ in $n$ variables which consists of
  polynomials from $B$ repeated $n - 1$ times with $x_2$ replaced by
  unchanged in the first copy, replaced by $x_3$ in the second, by
  $x_4$ in the third, $\cdots$, by $x_n$ in the last,
  respectively. The isolated solutions of system $C_n$ form a
  $n$-dimensional lattice consisting of $2(k-1)^{n-1}d^n$ points (there
  are $2d$ series each with a fixed value of coordinate $x_1$
  containing $((k- 1)d)^{n-1}$ isolated points).

  \begin{remark}
    System $A$ in 2 variables consisting of $k$ cubic tropical
    polynomials, has a linear in $k$ number of solutions. On the
    contrary, one can prove that a system in 2 variables consisting of
    an arbitrary number of quadratic tropical polynomials, has at most
    $72$ solutions.
  \end{remark}
\end{proof}
\end{theorem}
\section{Compatification of Tropical Prevarieties}
\label{section-compactification}
In this section we will show that for any tropical prevariety $V$ we
can build a compact tropical prevariety being homotopy equivalent to
$V$. This technique can be used to reduce the problem of estimating the
number of connected components of tropical prevarieties to the case of
compact prevarieties.

We will use the following theorem:

\begin{theorem}
  \label{theorem-bounding-cube}
  Given a tropical prevariety $V$ we can find a constant $s$ such that
  the intersection of $V$ and a cube with the side equal to $s$ and
  centered at the origin would be homotopy equivalent to  $V$.
  In this theorem we allow the prevariety to be given by a system with
  infinite coefficients.
\end{theorem}

This theorem can be viewed as a simplification of
Lemma~9 proved in \cite{grigoriev1988solving}, or it can be
proved directly with the help of Cramer's rule in the same way as it
was used in Theorem~\ref{stability-criteria}.

\begin{theorem}
  \label{theorem-compactification}
  Consider a tropical prevariety $V$ given by a tropical system $A$ in
  $n$ variables. We can add $2n$ extra variables and $4n$ extra
  polynomials which being added to system $A$ will form a system $B$
  which determines a compact tropical prevariety $W$ being a homotopy
  equivalent to $V$.

  \begin{proof}
    Let $2s$ be a side of the cube from
    Theorem~\ref{theorem-bounding-cube}.  For each variable $x_i$ we
    will add two variables: $u_i$ and $v_i$; and four (linear)
    polynomials: $x_i \oplus u_i$, $x_i \oplus u_i \oplus s$, $v_i
    \oplus -s$ and $x_i \oplus v_i \oplus -s$.  The first two
    polynomials will guarantee that $u_i = x_i \le s$, and the last
    two will guarantee that $x_i \ge v_i = -s$. Therefore, $W$ is
    homeomorphic to $V\cap [-s,s]^n$. Prevariety $W$ is compact and by
    Theorem~\ref{theorem-bounding-cube} it is homotopy equivalent to
    $V$. While there are infinite coefficients in the added
    polynomials this is not a problem as we can assume that all
    infinite coefficients are equal to $M + 2s^d$, where $M$ is a
    maximal coefficient occurring in the polynomials from $A$ and $d$
    is a maximal degree of polynomials from $A$. In that case monomial
    with this coefficient could not be minimal as all variables are
    guaranteed to be less or equal than $s$.
  \end{proof}
\end{theorem}

All results from the previous sections required all the coefficients
in a tropical system to be finite (as they were based on Bezout's
theorem, which was proved only in the case of finite
coefficients). However Theorem~\ref{theorem-compactification} gives us
the following generalizations to the case of infinite coefficients:

\begin{corollary}
  \label{corollary-connected-bound-infinite}
  The number of connected components of a tropical prevariety given by
  $k$ polynomials with tropical degrees bounded by $d$ and with allowed infinite coefficients is not greater
  than $\frac{d^{3n}}{k + n + 1} {{k + 7n - 1} \choose 3n}$.
\end{corollary}

\begin{corollary}\label{isolated}
  The number of isolated solutions of an overdetermined tropical
  polynomial system of $k$ equations in $n$ variables with allowed
  infinite coefficients  is not greater than $\frac{{k + 4n
      \choose 3n}}{(k + n + 1)}D$, where $D$ is the product of $n$
  greatest degrees of the given polynomials.
\end{corollary}

%\newpage
 {\bf Acknowledgements}. The authors are grateful to the Max-Planck
Institut f\"ur Mathematik, Bonn for its hospitality during writing
this paper.

\end{document}